# An Anarchist Approach to the Undergraduate Mathematics Curriculum

*Vincent Bouchard[1] and Asia Matthews[2]*

## Abstract

Contemporary anarchism centers around three tenets: (1) a constant challenge of and resistance to all forms of domination, (2) so-called "prefigurative politics", in which all decisions are made in a manner that is consistent with a set of non-hierarchical values such as equality, decentralization and voluntary cooperation, (3) a focus on diversity and open-endedness (Gordon, 2008). Within this philosophy the notion of end goals becomes moot; progress, then, is measured by process, in which the values of diversity, pluralism, cooperation, autonomy and experimentation are celebrated. In this perspective piece we propose anarchism as a philosophical framework to address the perceived cognitive dissonances of the current undergraduate mathematics curriculum. Are learning outcomes appropriate in an anarchist approach to education? How can we address the power dynamics of grading and assessment? How can assessment be done in the context of a process-based and horizontal approach that celebrates diversity and autonomy? Should grades be used, and if so, how could they be assigned non-hierarchically? At its core, anarchism aims at aligning thoughts and actions, and we argue that an anarchist viewpoint on undergraduate mathematics addresses the cognitive dissonances that currently plague our curriculum. We propose food for thought for individual instructors' practice, including ideas for incremental and large-scale changes.

___________________

[1]*Department of Mathematical and Statistical Sciences, University of Alberta, Edmonton, Canada*

Email: vincent.bouchard@ualberta.ca

[2]*Quest University, Squamish, Canada* and *Department of Mathematics and Statistics, University of Victoria, Victoria, Canada*

Email: dr.asia.matthews@gmail.com

# 1. Introduction: A Cognitively Dissonant Curriculum

You are trying to moderate your alcohol consumption. You're invited to a social event. You have a drink. Your friend offers a second round. You say yes, but you know that you shouldn't have a second drink, because, well, you're trying to moderate your consumption. You feel discomfort and guilt. This is cognitive dissonance.

Cognitive dissonance can be described as a "mental conflict that occurs when your beliefs don't line up with your actions".[3] It is perhaps one of the most uncomfortable states of mind. To avoid the resulting internal conflict, we often try to rationalize our decisions and come up with bogus explanations. But it never really works; the deep feeling of unease remains, which makes the whole process counterproductive. Why not instead simply align our actions with our values?

The authors, as many mathematicians and undergraduate educators, have experienced cognitive dissonance in the way that mathematics is taught at the undergraduate level. This is a common lament: see, for example, Lockhart (2002). We present here a collection of expressions that highlight the unease that we have experienced.

We want our students to learn mathematical skills, such as abstraction, problem-solving, creativity, rigour, logical thinking, collaboration, etc. Yet most of us focus our teaching and assessment almost exclusively on mathematical content. We do not teach students how to do math; we teach them how math was done by others. We teach the mathematical product instead of the mathematical process.

We proclaim great importance to equity, diversity and inclusion. Yet we teach math in a very hierarchical way, as if knowledge trickles down from the mathematician to "his" pupils. And, as the field has been dominated by males, and typically white males in North America and Europe, this knowledge is filtered, molded and aimed in ways which form barriers for those not identified in these dominating categories.

We believe that the mathematical process that leads to an answer is as important, if not more, than the answer itself. Yet we design our courses by focusing on learning outcomes rather than processes.

We value mathematics as a way of thinking and approaching the world; we would like to overcome the deeply ingrained fear of mathematics in our society. Yet our courses often act as "gatekeepers", with low averages and high failure rates. Students in non-mathematics degrees often see mathematics courses as tedious barriers to completion of their degree.

We experience strong emotions when we see a beautiful mathematical argument. We hope to share our excitement with our students. Yet we teach mathematical content mostly as a dry, heart-disconnected subject.

---

[3] https://www.psycom.net/cognitive-dissonance

We see mathematics as a creative endeavour, yet we mostly teach algorithmic procedures. We believe in the "throw ideas at the wall and see what sticks" creative process, yet we teach content that has been manicured to make it as clean and structurally sound as possible.

And we could go on and on. The cognitive dissonance of our undergraduate mathematics curriculum is all-pervasive and so ubiquitous that we barely see it. We rationalize it using all kinds of arguments; the need to prepare students for the job market, the desire to have clear "outcomes" that can be easily assessed, and so on and so forth. Those are all valid points, but ultimately, the internal conflict remains, and we feel uneasy and uncomfortable. This is cognitive dissonance.

This issue has certainly been discussed in the literature. Long-standing enculturation and prejudice affect the choices teachers make in the classroom (Reinholz, 2023), and though many have sought to understand the complex psychological relationships in the mathematics classroom (Skemp, 1971), continuous personal, classroom environment, institutional, cultural, sociological, and global changes make understanding elusive and solutions narrowly applicable. Feeding the classroom is an aging mathematics curriculum, designed by humans holding their own beliefs and values which may or may not be shared by each individual mathematician and educator in the classroom. Furthermore, the complex relationship between student-centered educational values, socially- and globally-focused institutional values, and the mathematics curriculum, creates situations in which educators must act in ways which are at odds with their values (Chronaki & Yolcu, 2021).

In this perspective piece we propose anarchism as a philosophical framework to address the perceived cognitive dissonance of the current undergraduate mathematics curriculum. Anarchism is a philosophy that, at its core, aims at resolving cognitive dissonances by aligning thoughts and actions. Furthermore, we argue that mathematics, as a practice, is inherently anarchist; thus, if our teaching followed anarchist principles, cognitive dissonance may be reduced. Ultimately, we ask the simple yet profound question: what if we aligned our undergraduate mathematics curriculum with our values? Through the lens of anarchism we suggest ways in which instructors might reconceive the curriculum and the mathematics classroom to bring them more in line with their values and minimize the experience of cognitive dissonance.

## 2. Anarchism and Mathematics

### 2.1. Anarchism

Anarchist ideas have existed in human society for a very long time. Anarchism as a political philosophy developed in the late 18th and 19th centuries through the writings and actions of Godwin, Proudhon, Bakunin, Kropotkin, and many others (Rocker, 1989). Many different anarchist schools of thought have emerged over time, from very individualistic trends to deeply socialistic approaches, so it is of crucial importance to specify what we mean by anarchism here, especially given that anarchism is so often misconstrued in society.

We focus on contemporary anarchism, which centers around three tenets: (1) a constant challenge of and resistance to all forms of domination, (2) so-called "prefigurative politics", in which all decisions are made in a manner that is consistent with a set of non-hierarchical values such as equality, decentralization and voluntary cooperation, and (3) a focus on diversity and open-endedness (Gordon, 2008). Within this philosophy, the notion of end goals becomes moot; progress, then, is measured by process, in which the values of diversity, pluralism, cooperation, autonomy and experimentation are celebrated.

The idea of prefigurative politics, which is core to anarchism, is the fundamental remedy to cognitive dissonance. It is the concept that we should act, behave and make decisions that are always consistent with what we value and hope to achieve. As Emma Goldman (1923) wrote:

> There is no greater fallacy than the belief that aims and purposes are one thing, while methods and tactics are another. This conception is a potent menace to social regeneration. All human experience teaches that methods and means cannot be separated from the ultimate aim. The means employed become, through individual habit and social practice, part and parcel of the final purpose; they influence it, modify it, and presently the aims and means become identical… the ethical values which the revolution is to establish in the new society must be initiated with the revolutionary activities of the so-called transitional period. The latter can serve as a real and dependable bridge to the better life only if built of the same material as the life to be achieved. (Afterword)

Prefigurative politics is the opposite of the Machiavellian saying that "the end justifies the means". For plants to grow, seeds must be planted today. It is the "this is what democracy looks like" chant at demonstrations and direct actions; it is the famous prediction by Bakunin that if Marxist revolutionaries gained power, they would not create an equal, free, communist society, but would rather end up new tyrants of the working class. We suggest that it is also what may be applied to address the cognitive dissonance of our undergraduate mathematics curriculum.

The two other central tenets of contemporary anarchism are also keys to rebuilding the undergraduate mathematics curriculum. Indeed, as we now explain, the fundamental challenge of all forms of domination, and the focus on diversity and open-endedness, play a prominent role in the way that mathematics is done, and in how research in mathematics is conducted, and as such, should also be reflected in the way mathematics is taught.

*2.2. Mathematics as an Anarchist Endeavour*

Although some aspects of mathematical practice and education are permeated with hierarchical dominance and unethical constraints, we argue that mathematics, as a disciplined pursuit, is inherently anarchist: it is open-ended and progresses through a plurality of approaches, while deductive reasoning is the ultimate authority.

Mathematics is open-ended. As Cantor (1883) remarked, "the very essence of mathematics is its freedom," and open-endedness is the foundation of mathematical research. When starting a research project, the mathematician does not have predetermined outcomes; instead, they ask a question and try to provide answers, remaining true to our core value of clear, logical reasoning every step of the way (prefigurative politics). It is an organic process which often leads to places that could not have been imagined.

Mathematics is diverse. The plurality of approaches is also crucial, as often revisiting a mathematical concept from a different perspective leads to deep realizations and new connections. Mathematics in its purest form is a celebration of unity in diversity, with open-endedness being the foundation of mathematical research.

Mathematics is free to all. The publication process in mathematics shares some anarchist values in that written mathematics, reviewed by subject experts of logic and rigour, is shared quite freely within a smaller community of experts and more widely through the open-source platform arXiv, which is freely available.

Mathematics rejects power structures. The mathematical process relies solely on deductive reasoning: what is true mathematically should be decided by logical reasoning, not appeal to authority. This is true both in mathematical practice, and should also be the case in mathematics education.

> The final criterion of any piece of mathematics is consistency: with itself, that is to say internal consistency, or with the larger mathematical system of which it forms part. Whether this consistency exists is a matter for agreement between one mathematician and another, and between teacher and learner. [...] What is more, the criterion is implicitly accepted as binding by teachers and students alike. [...] In mathematics, perhaps more than any other subject, the learning process depends on agreement, and this agreement rests on pure reason. (Skemp, 1971, p. 116)

Acceptance of mathematical truth is based on consensus within the mathematical community.[4] Whether or not this practice is reflected in mathematics education, mathematics is done by humans, and communities of practitioners and learners bring their own biases and inequities. Reinholz (2023) cautions us to pay close attention to the interwoven nature of mathematics and the practice of mathematics, remarking that a perception of mathematics as "rational, perfect, objective, and complete" is simplistic at best. Mathematics, they continue, "has a special status in society, which uniquely positions it to signify – and thus perpetuate – inequity." While mathematics itself may present objective truth, the discipline, application, and education of mathematics, "is still imperfect, and influenced by human judgments, flaws, and biases" (p. 26).

---

[4] It will be interesting to see how this is changing with the development of computer-based formal proof assistants such as Lean, a functional programming language and interactive theorem prover (https://lean-lang.org/).

Nevertheless, mathematics, mathematicians, and the mathematical research process share values and practices of openness, fairness, and logical reasoning that are central to anarchism. Therefore we argue that the three central tenets of anarchism may provide a cure to the diseases that currently plague undergraduate mathematics education.

*2.3. Anarchism, Chaos, and Expertise*

We now consider mathematics education, raising and addressing a number of questions. If the anarchist tenets of decentralization and voluntary cooperation support individual pursuit of knowledge acquisition, does anarchism prohibit an organized and predetermined education system? Does an anarchist perspective on undergraduate mathematics education preclude even having systemic goals? More extremely, is it anti-intellectual, going in the direction of a "cult of ignorance"?

First, the anarchist perspective is not the absence of order. It is not chaos, even though it is often portrayed as such in popular culture. The authors do not propose that undergraduate mathematics education become a non-system in which there is no curriculum, no goals, and no direction. We simply propose aligning the curriculum (design and execution) with the three central tenets of contemporary anarchism, a focus on process, open-endedness and diversity, and resistance to hierarchies.

But does challenging hierarchies imply that a student's opinion on a particular mathematical topic is always just as valid as the teacher's opinion? Certainly not. Challenging hierarchies does not mean the "death of expertise" (Nichols, 2017). But this is a sensitive issue, particularly in this day and time, so let us expand a bit more.

An anarchist approach may see the teacher as a guide, whose expertise is available for students to draw upon. This is not the same as being an authority, whose word must be taken for granted. This can be a difficult line to walk, especially in the elementary and secondary mathematics classrooms, where two kinds of authority – expertise and classroom discipline – can be at odds (Skemp, 1971, Chapter 7). In this dual role, it is the responsibility of the teacher to show how mathematical truth is determined by reason not hierarchy. Thus, no word should be taken for granted, everything should in principle be questioned and proved. Furthermore, mathematical expertise is a guide for mathematical processes such as determining validity and developing intuition for mathematical courses of action.

The critical distinction to make is between expertise and power. In mathematics, the separation could not be clearer: the expert mathematician will be the first to emphasize the importance of scrutinizing mathematical statements, regardless of who proposed the claim, and making sure that correct proofs exist. (Think of Fermat's last theorem…) No mathematician would ever accept an appeal to authority as a valid mathematical proof.

Ultimately, the anarchist approach provides a clear and practical distinction: expertise is essential and respected; power should be challenged. It may be illuminating to think of the relation between a patient and a physician. The physician has expertise, and after diagnosis may recommend a course of action to the patient. Most patients would then agree and follow the recommendation, respecting the expertise of the physician. However, ultimately the doctor does not have the power to impose the course of action on the patient (assuming the patient is conscious and coherent). There is a clear distinction between expertise and power. An anarchist approach to mathematics education should also make such a distinction, with the expert making recommendations to students to help them flourish and become experts in turn.

## 3. Refocusing the Undergraduate Mathematics Curriculum

What would an anarchist approach to undergraduate mathematics education look like?

First and foremost, we remind the reader that there is no anarchist blueprint of what an anarchist society would look like, and there will never be.

> Anarchism is no patent solution for all human problems, no Utopia of a perfect social order, as it has so often been called, since on principle it rejects all absolute schemes and concepts. It does not believe in any absolute truth, or in definite final goals for human development, but in an unlimited perfectibility of social arrangements and human living conditions, which are always straining after higher forms of expression, and to which for this reason one can assign no definite terminus nor set any fixed goal. (Rocker, 1989, p.30)

An anarchist blueprint would be implicitly hierarchical, as if the author of the blueprint knows better than everyone else. It would also run against the idea of prefigurative politics, by first proposing the goal instead of focusing on the process. As Chomsky (1987) wrote, this makes anarchism "an unending struggle, since progress in achieving a more just society will lead to new insight and understanding of forms of oppression that may be concealed in traditional practice and consciousness" (p. 50).

Instead, we propose a few avenues to inspire deeper thought by investigating how the three central tenets of modern anarchism could be applied to help resolve some of the fundamental dissonances that we identified earlier. Specifically, the anarchical principles which resist domination, hold non-hierarchical values, and esteem diversity and open-endedness, suggest a change in the focus of the undergraduate curriculum from content to skills, from outcomes to process, and from hierarchical to horizontal teaching structures. As progress is made towards achieving a curriculum which aligns more closely to values held by contemporary mathematicians and institutions, experiences will lead to new insights, new ways of understanding ingrained forms of domination, new questions, and new paths forward. We hope that such a focus will provide those who experience it a lightening of the burden of cognitive dissonance.

### 3.1. Skills Over Content

*We want our students to learn mathematical skills, such as abstraction, problem-solving, creativity, rigour, logical thinking, collaboration, etc. Yet most of us focus our teaching and assessment almost exclusively on mathematical content.*

The importance that contemporary mathematics curriculum places on knowledge can be seen in the way that undergraduate courses are named: by content category, e.g. Linear Algebra, Calculus. However, *doing* mathematics and *learning* to do mathematics depends on cognitive and metacognitive *skills* such as creative thinking and self-regulation, skills which are rarely explicit in mathematics curriculum.

The OECD Learning Compass 2030 (2018) distinguishes between three different types of skills:

- Cognitive and meta-cognitive skills, which include critical thinking, creative thinking, learning-to-learn and self-regulation;
- Social and emotional skills, which include empathy, self-efficacy, responsibility and collaboration;
- Practical and physical skills, which include using new information and communication technology devices.

According to the OECD (ibid.), skills are part of a holistic concept of competency, involving knowledge, skills, attitudes and values which are developed interdependently. Competency, we suggest, may be understood as an essential step toward mastery in mathematics which has been conceptualized in many ways. Various countries throughout recent history, Niss, et al., (2017) explain, have seen shifting emphasis in curricular focus on: knowledge, skills and techniques, the process of doing mathematics, and mathematical thinking. None of these, they explain, exist independently, and this is what is intended by the idea of the more holistic approach of competencies. From the mid-1980s through the early 2000s we saw the instrumental Common Core State Standards for Mathematics (CCSSM) highlight competencies in the school mathematics curriculum, though it remains unclear if these changes have infiltrated the undergraduate mathematics curriculum.

Yet, even if the *official* curriculum emphasizes competencies or elements beyond knowledge, curricular intentions can transform as they move into the *enacted* curriculum (Remillard and Heck, 2014). While the authors are not aware of supporting research, our experiences as members of the collection of undergraduate faculty in both the United States and Canada have led us both to believe that the majority of *enacted* undergraduate curriculum is focused primarily on knowledge and understanding of content (definitions, concepts, theorems) and on basic cognitive skills such as algorithmic procedures and techniques. It is not all doom and gloom; some faculty do emphasize metacognitive skills, and practical and physical skills are becoming more important with technological advances. Collaboration seems to be gaining traction as well. Social and emotional skills however have little to no explicit curricular emphasis in undergraduate mathematics despite contemporary mathematics education research showing strong connections between safe communities and the ability to learn (see, for example, Chapter 4 of Reinholz, 2023).

Ultimately, many mathematicians, educators, and educational institutions greatly value cognitive and metacognitive skills such as abstraction, creativity, problem-solving, collaboration, communication, etc. Francis Su (2020) refers to "virtues" as being built by engaging in exploration, meaning, play, beauty, permanence, truth, struggle, power, justice, freedom, community, love. In a recent interview (Gu & Bishop, 2022) Su shares, "Virtues make your life richer. These are the things that you'll carry with you the rest of your life, no matter what profession you go into, and no matter where your life takes you." Other stakeholders, such as employers, also agree: virtues developed through mathematics education and mathematical reasoning skills can help students thrive in our world and shape a better future (OECD, 2018).

With an emphasis on prefigurative politics, an anarchist approach to mathematics education would suggest that if we do really value cognitive and metacognitive, social and emotional, and practical and physical skills, then at every step of the way we should orient our teaching and curriculum to focus on these skills, instead of viewing them as byproducts of a curriculum designed around content. From an anarchist perspective, the process of doing is paramount to that which is produced; this reflects an assertion by Skemp (1971) that teaching should focus more on the process of mathematical thinking and less on the product of mathematical thought (pp.13-14).

As an example, consider what a difference it would make if our courses were organized around valuable mathematical thinking processes instead of content. Rather than classifying our courses by products of mathematical thought, courses could be offered by processes of mathematical thinking. An undergraduate curriculum might then be structured around courses such as:

- Abstracting
- Communicating and collaborating
- Computing
- Acting ethically and socially responsible
- Problem-solving
- Writing proof and using formalism
- Researching

This is, in some sense, a liberal education within mathematics focused on building strong foundations of ways of thinking about mathematics, using the traditional mathematical content as the medium through which to explore these processes. The instructor's goal would be to help students become comfortable and familiar with applying these processes, not only to construct mathematical arguments, but in everyday life as well. Core mathematical knowledge and skills would still be covered in such courses (such as the notion of derivatives etc.), but the knowledge would not be taught solely for the sake of knowledge; rather, it would be used as context to put in practice the appropriate skills.

This kind of sweeping change is unlikely to occur quickly, and the anarchist perspective can provide some direction on a smaller scale as well. Changes in individual courses can also transform the focus from knowledge and basic cognitive skills to include metacognitive, social and emotional skills, or can choose themes from the list above to reinforce throughout a course. For example, regular tasks which require students to develop visual representations of mathematical ideas can be aimed at abstracting and communicating. Frequent group work and peer feedback can be designed well to probe deeply into students' metacognitive process as well as their collaboration and social responsibility. Regular tasks which require students to extend provided problems to create new, related problems (see the literature on problem posing) is one way to help them develop their creativity and provide a taste of a beginning to the research process.

### *3.2. Horizontal Teaching Rather Than Hierarchy*

*We proclaim great importance to equity, diversity and inclusion. Yet we teach math in a very hierarchical way, as if knowledge trickles down from the mathematician to "his" pupils.*

The hierarchical structure of undergraduate education takes many forms. In the curriculum, for example, the mathematician, the subject expert, specifies the precise path that the learner should take. The instructor is also responsible for evaluating students and assigning grades, thus reinforcing the power structure. In the teacher-class relationship, the role of teacher-as-manager influences the individuals who get to participate in classroom activities and discussions. An instructor's choices are grounded in bias, often leading to continued marginalization, inequity, and a reinforced hierarchical classroom structure.

> In addition to biases, Discourses position some students as more capable at mathematics—and in the absence of an instructor intentionally mediating the discussion—these students tend to dominate classroom discussions. For example, masculine Discourses in mathematics produce a classroom context in which men exhibit louder, aggressive, and competitive behaviors, which in turn, means that they take up a disproportionate amount of talk time and learning opportunities. (Reinholz, 2023, p. 29)

Hierarchies perpetuate not only in the teacher-student relationship but also among students. As Reinholz (2023) explains, hierarchies are grounded in Discourses, "a collection of symbols, signs, artifacts, and other cultural representations that work together to constitute the social world." Social Discourses separate groups into categories and then associate ranking, and thus a power imbalance, to these groups. The categories can be, for example, racist ("Asians are better at math"), sexist ("women are less logical"), or ableist ("students with accommodations are deficient"). "[T]hese Discourses," Reinholz continues, "locate problems within students and their communities, while obscuring the role of systems of oppression that perpetuate inequity. [...] [M]athematics is not devoid of culture, but rather, is racialized, gendered, and so forth" (p. 28).

Mathematics educators have an opportunity – a responsibility – to address inequity and its implicitly associated hierarchies. From an anarchist perspective, a key to addressing this is to collectively (students and faculty together) work to create an open, inclusive, respectful and non-hierarchical learning environment which more closely mirrors the inherently non-hierarchical mathematical process itself. Curricular changes are necessary, from the anarchist perspective, as well. An idealistic view might suggest that students and teachers agree on a curriculum together at the beginning of every course. Practically, this is unlikely to happen within the confines of the current system. But this is exactly what the anarchist perspective provides - an understanding that all systems are flawed, and a way to take steps forward.

The non-hierarchical nature of the mathematics subject can be a guide for addressing social inequity and educating global citizens who recognize and address irrational power imbalances outside the mathematics classroom. As argued before, there is no denying that expertise and established knowledge should be respected; all opinions are not equally valuable in mathematics. But we make a distinction between expertise and power. The expert mathematician knows more about their subject than the first-year student, but this does not need to translate into a power dynamic. As subject experts, we can make choices which provide students guidance and help them develop into independent critical thinkers, not followers.

> The student has no need to accept anything which is not agreeable to his own intelligence – ideally he has a duty not to. And it is by the exercise of the teacher's intelligence, and not by his prestige, eloquence, or tyranny, that the students should be led to agree with him. The teaching and learning of mathematics should thus be an interaction between intelligences, each respecting that of the other. (Skemp, 1971, pp. 116-117)

Instructors can make many curricular and pedagogical choices to emphasize horizontal structures in their teaching. One may choose to address colonization, while another might focus on power imbalance in classroom discussion. Others may make more drastic changes, for example from lecture-based to inquiry-based pedagogies. But what of changes in the curriculum?

One obvious hierarchical structure in undergraduate education is in grade assignment – at the end of a course, the teacher, the subject expert, makes a judgment about the "level" the student is at. The process is implicitly inequitable and hierarchical. As Reinholz comments, "grading is primarily a mechanism to create and sustain hierarchies. [...] Because mathematics is seen as a signifier of intelligence, [...] testing regimes [...] sort students and constrain future opportunities" (Reinholz, 2023, p. 89). A detailed look into assessment and hierarchies in mathematics education can be found in Chapter 6 of Reinholz (2023).

"Ungrading" is an attempt at addressing this problematic issue:

> The word "ungrading" means raising an eyebrow at grades as a systemic practice, distinct from simply "not grading". The word is a present participle, an ongoing process, not a static set of practices. Ungrading is a systemic critique, a series of conversations we have about grades, ideally drawing students into those conversations with the goal of engaging them as full agents in their own education. For me, there aren't a discrete set of best practices for ungrading, because different students learn in different ways at different times with different teachers in different disciplines at different institutions. […] Grades are inequitable. As they are increasingly centered at our institutions and within our educational technologies (like the learning management system), the inequities of grades are exacerbated, and our most marginalized students are further marginalized. (Stommel, 2023, p. 6)

Sounds familiar? Ungrading, which puts an emphasis on process, a celebration of diversity, and a challenge of hierarchies, fits well within the main themes of anarchism. In most universities, a letter grade is required at the end of a course. From both the ungrading and the anarchist perspectives, we, as instructors and institutions, begin by questioning if the way that we are assigning grades agrees with our values. This approach is not new, but is not very common in mathematics courses. Perhaps it is because many educators are worried about "quality control," given the scaffolding of courses in mathematics. Perhaps it is our proclivity toward numbers.

Some undergraduate instructors have freedom in their choice of assessment while others are constrained to departmental curricular choices. For the first, there are many developing assessment processes that align more closely with mathematical and institutional values, including self-assessment, mastery grading, holistic grading, and non-numeric grading. Faculty working within a more rigid departmental assessment structure, can develop hybrid systems. For example, a hybrid approach to ungrading and self-assessment may involve having students self-assign their own letter grade for term work evidenced by a student portfolio, with the instructor assigning a separate letter grade based on performance on a final assessment (which could be an oral exam, a written final exam, …). The final letter grade could then be the average of the two, with the instructor's grade superseding the student's grade if it is higher (so that students who are too harsh with themselves are not penalized). In effect, this amounts to implementing a system in which letter grades are assigned collaboratively between students and instructors, still ensuring a level of quality control for further courses.[5]

---

[5] The first author (VB) successfully implemented such a system in a large (150 students) undergraduate mathematics course. In a large classroom setting, it is challenging to provide adequate feedback throughout the semester so that students can accurately self-assess their own learning process. Nevertheless, the success of the experiment shows that it is certainly possible!

Or perhaps we should get rid of grades altogether?[6] The undergraduate mathematics curriculum attempts to provide clear learning outcomes which are intended to be reached by the end of a course. Grades implicitly assume measuring a student's output in comparison to this predefined list of outcomes, so from an anarchist perspective grading goes against the tenet of open-endedness. This brings us to the next point: does an anarchist perspective preclude a focus on outcomes in course design?

### *3.3. Mathematical Thinking as Outcomes*

*We believe that the mathematical process that leads to an answer is as important, if not more important, than the answer itself. Yet we design our courses by focusing on learning outcomes rather than processes.*

Mathematical thinking processes – analytic thinking, rigour, precision, argumentation, the proper application of procedures, conjecturing, creativity – are crucially important and fundamental to doing mathematics and to solving problems with mathematics. It is of course also critical to have knowledge of concepts and procedures (content) to apply in situations. These are both essential and interwoven parts of mathematics. The current undergraduate curriculum, however, values mathematical knowledge over the process of thinking mathematically. This is evidenced by the pre-determination of learning outcomes, by standardized assessment practices, and by the titles of undergraduate courses themselves. As countless articles and editorials lament, the curriculum is not serving mathematics education stakeholders well.

Mathematics is an open-ended process, celebrating the diversity and plurality of approaches. However, in contemporary mathematics education, the content of a course is identified by learning outcomes – measurable statements that describe the knowledge or skills that students should acquire by the end of the course. This practice is a focus on the endpoint, the goal to be reached, with little emphasis on the process towards this goal. It is no wonder that many students perceive mathematics as a set of rules and procedures intended on reaching established predefined outcomes.

Learning outcomes go against anarchist principles, are fundamentally goal-oriented, and reward uniformity. A particular teacher could teach a course with sole emphasis to get their students to perform well on the examinations, which are designed to assess the outcomes of the course, with little regard to whether students are becoming better at actually doing

---

[6] Grading has been shown to be detrimental to the learning process. However there is more to the story: "Indeed, one could argue that higher grading standards may motivate some students to work harder to achieve high grades and at the same time serve as a disincentive to other students who may now view certain grades as out of reach. (…) The effect of grading standard on performance is theoretically and often empirically ambiguous. Despite this uncertainty, there is little doubt that incentives serve to motivate behavior" (Elikai & Schuhmann, 2010).

mathematics. Teaching for the exam is the Machiavellian "the learning outcomes justify the means".[7]

Furthermore, learning outcomes are dominantly imposed by centralized authority, with no input from students. This goes against the anarchist principle of decentralized authority. In fact, many university mathematics professors hate learning outcomes, do not value them, and write them poorly. Why is this? Perhaps because of an experience of cognitive dissonance between valuing the processes of thinking mathematically and the practice of goals aimed at the product of mathematical thought.

From an anarchist perspective, learning outcomes are hierarchical and closed-ended because they are designed and imposed on students from an authority and focus on goals rather than the process of learning. Thus, it appears that from this perspective we must discard all learning outcomes. On the other hand, it is often said that learning outcomes are helpful for students, as they provide a clear and expected learning path. How shall we proceed?

By definition, a course has a beginning and an end. Is it possible to structure a course while adhering to the anarchical tenet of open-endedness? What could an "open-ended course" mean? Could courses be designed with flexible learning outcomes? Perhaps with a list of outcomes that students can select from?

Let us first suppose that we keep learning outcomes as the main focus of course design, as that is indeed the status quo. One value-driven choice could be to re-envision an "outcome" as a suggested path toward learning how to think mathematically rather than a collection of knowledge and skills to acquire. To implement some level of open-endedness, and reduce the centralized authority implicit in designing learning outcomes, we could allow for the possibility that, as a student progresses through a course, they may want to change their goals. Or perhaps we do not specify learning outcomes for a course but instead ask students to reflect (a metacognitive skill) and write the experienced learning outcomes for themselves at the end of the course. It is, after all, at the end of a course that one can reflect on the experience as a whole and describe and evaluate the learning that took place.

A focus on mathematical thinking processes over outcomes can also lead to changes on a smaller scale. We can design single courses that focus on the process of mathematics instead of the mathematical product. For example, a first-year modeling course following the Society for Industrial and Applied Mathematics Math Modelling handbooks (SIAM, 2014, 2018) focuses on modelling processes (e.g. Defining the Problem Statement, Making Assumptions, Communicating Results) over any pre-determined procedural techniques[8].

---

[7] "Teaching to the test" has been widely discussed in mathematics education. See for instance (Klymchuk & Sangwin, 2020).

[8] The second author (AM) has done just this and observed exceptional growth in students' critical thinking and handling of quantitative information. Students reflected that they became less afraid and more capable of doing mathematics. This course was inspired by Yvan Saint Aubin's modelling course offered at Université de Montréal.

Curricular changes on an even smaller scale might use various pedagogical methods within a course to emphasize open-endedness. An example is dialogue, which is not used often as a pedagogical tool in mathematics courses.

> Dialogue has no pre-defined direction, and the results can never be predicted in any absolute sense. (...) Thus, when talking about dialogues, we are ideally speaking of equality, complementarity and mutuality in the context of a communicative project but not of complete symmetry or overlap in perspectives or opinions. In an ideal dialogue there should be no use of overt power or force, no up-front persuasion of the other, and no-one should win at the expense of the other. The purpose or outcome of a dialogue should not be defined or decided on through authority. (Skovomose & Säljö, 2008, p. 36)

Or perhaps we discard outcomes altogether. A sweeping change could be to design courses to mimic the way that we do research in mathematics, which, as we saw, shares many anarchist values. In this way, courses could be designed around a set of questions –- questions which are open-ended, not questions that we already know the answer to – rather than outcomes. Questions such as "Is population growth an indicator for inequity?", or "What could go wrong with writing a mathematical argument in informal terms?" could aim at a variety of mathematical thinking processes as well as cognitive, metacognitive, social and emotional skills. A more decentralized choice could be structured around, say, 6 or 7 questions and students could choose 2 or 3 to focus their learning experience on.

As an example, a course built around questions rather than learning outcomes could then be focused on the mathematical thinking required to propose answers to the selected questions. Students would build a portfolio detailing their thought processes to answer the questions. The instructor's role would be to guide students in their inquiry and to advise on mathematical knowledge that could be used to propose deeper answers to the questions. If the questions for a given course are designed carefully, some of the relevant mathematical knowledge could be common to the questions. The instructor could then teach regular "modules" to students, covering mathematical knowledge that may be relevant for their investigations, but let the students figure out how they may use it to answer their own questions.

This goes in the direction of Freire's critical pedagogy and libertarian education (2005), in which the classroom should be a space for asking questions:

> Education must begin with the solution of the teacher-student contradiction, by reconciling the poles of the contradiction so that both are simultaneously teachers and students. [...] Liberating education consists in acts of cognition, not transferrals of information. It is a learning situation in which the cognizable object (far from being the end of the cognitive act) intermediates the cognitive actors – teacher on the one hand and students on the other. Accordingly, the practice of problem-posing education entails at the outset that the teacher-student contradiction to be resolved. Dialogical relations – indispensable to the capacity of cognitive actors to cooperate in perceiving the same cognizable object – are otherwise impossible. (pp. 72-80)

What does an anarchist perspective provide with respect to assessment without learning goals? A common argument in support of learning outcomes is that they provide clear, quantifiable goals that can be assessed. It is usually said that summative assessments should focus on measuring whether students have met the predefined learning outcomes.

> Narrow conceptions of mathematics content are often paired with narrow conceptions of understanding, as an answer either being "right or wrong" for a standardized assessment. These conceptions focus little on processes, and don't reflect authentic mathematics. When we move beyond covering content, we can attend to the types of mathematical practices that our students engage with. (Reinholz, 2023, p. 54)

Designing courses around processes or questions instead of goals begs a new type of assessment. Measuring students' level in comparison to learning outcomes predetermined by the instructor does not agree with anarchist ideals. If we adhere to these values, we need to be creative and find other ways of measuring students' progress. Fortunately, there are other approaches to assessing students' progress, beyond standard summative assessments. Again, frameworks such as mastery-based learning, or building a student portfolio, or even oral exams can be used to assess students' progress without having predefined learning outcomes. Peer review is also valuable for both parties and mirrors the authentic experience of research.

A curious byproduct of a curricular focus on processes instead of goals may be a new approach to the timeless issue of cheating. Cheating is inherently machiavellian; whatever the means, reach the goal. But if there is no goal, what would cheating even mean? If the whole course is based on the process, the means become indistinguishable from the aims. Whatever assessment framework would be used in a course that is designed around processes (or questions), as it would be based on measuring students' progress, the issue of cheating would presumably become much less prevalent.

Pushing the idea of open-endedness to its obvious conclusions suggests that we reconsider the idea of courses altogether. Perhaps undergraduate education should allow students to forge their own path? Or perhaps we do away with courses after the first or second year and students identify a research area or open-ended question that they are interested in, and then apply the skills that they have learned towards this question, under the supervision of subject experts. The learning environment would resemble the horizontal nature of research collaborations in mathematics; the interaction between the teacher and the learner would be one between intelligences, not between egos.

### 3.4. Affect in Mathematics

*We value mathematics as a way of thinking and approaching the world; we would like to overcome the deeply ingrained fear of mathematics in our society. Yet our courses often act as "gatekeepers", with low averages and high failure rates.*

What is the point of mathematics courses? It is really to make sure that students can take the derivative of an exceedingly complicated function? Or is it rather that they can apply mathematical processes such as problem-solving and abstraction to not only mathematical problems, but also all kinds of real-life situations?

The fear of "mathematics" is pervasive in our society. Students are scared of what they believe to be mathematics; high failure rates and low GPAs in mathematics courses also contribute to the fear of math. It would require an essay on its own to uncover the roots of this fear in our current society, but what is clear is that the current mathematics curriculum, with its emphasis on dry mathematical concepts with little space or time for joy and excitement, is counterproductive.

Let us consider the issue of failing. From the perspective of prefigurative politics, assessment of learning should include both instructors and students. Taking the viewpoint of designing courses based on processes or questions rather than outcomes, it may make sense not to fail students. By the end of a course students would necessarily have progressed at different paces, and reached different levels of understanding. Whatever system would be used to "grade" this progress, it would measure how far students have progressed. But from this perspective, no one would fail, as there is no predetermined goal, only different paces and learning progressions. Anarchism promotes unity in diversity, with an emphasis on open-endedness.

Importantly, what institutional and curricular problems does this cause in the long run? It is often said that the scaffolding of courses in mathematics education requires assessments based on learning outcomes to establish whether a student may proceed with another course which relies on a particular set of knowledge or skills. However, this may not be necessary. Perhaps a cooperative assessment process would provide the student with reasons that could help them make responsible choices about the courses they sign up for. Furthermore, many institutions are beginning to explore and promote personalized learning, and these types of curricular changes suggested by the anarchist viewpoint agree well with such values and support the ideas of diversity and open-endedness. Under such systematic changes, students who are not ready to proceed in mathematics will either not proceed, or will take longer to proceed.

*We experience strong emotions when we see a beautiful mathematical argument. We hope to share our excitement with our students. Yet we teach mathematical content mostly as a dry, heart-disconnected subject.*

Mathematics is an art. Mathematicians write novels, which include events they call theorems. Written mathematics is storytelling, but for a particular type of stories around patterns, symmetries, numbers, shapes, and in the particular language of logic. A novel is not great because of its impeccable grammar, though this plays a role; similarly, the beauty of a mathematical argument does not lie solely in that it is logically sound. It is difficult to pinpoint what makes a novel or a mathematical argument beautiful, but once experienced, the beauty is undeniable.

Beauty is evocative of strong emotions, and this is a phenomenon experienced and historically commented on by mathematicians. Frequent writing about aesthetic and intuition in mathematical activity abounds – see, for example, in Breitenbach and Rizza's (2018) introduction to the special issue on aesthetics in mathematics:

> Mathematicians often appreciate the beauty and elegance of particular theorems, proofs, and definitions, attaching importance not only to the truth but also to the aesthetic merit of their work. As Henri Poincaré [1930, p. 59] put it, mathematical beauty is a 'real aesthetic feeling that all true mathematicians recognise'. Others went further, regarding mathematical beauty as a key motivation driving the formulation of mathematical proofs and even as a criterion for choosing one proof over another. As Hermann Weyl famously and provocatively declared, 'My work always tried to unite the true with the beautiful, but when I had to choose one or the other, I usually chose the beautiful' (cited [Chandrasekhar, 1987, p. 52]).

In addition to the product of mathematical expression, the mathematical process itself can be very emotional. Many mathematicians comment on experiencing frustration when stuck in a research problem but a powerful ecstasy once a beautiful solution is found. Along the way, as we struggle, the little "aha" moments - the moments of sudden insight and discovery - keep us invested.

There has been a great deal of conversation in the mathematics education literature about affect in learning mathematics (e.g. McLeod, 1992), some of which points to a need to focus curriculum on more than the purely cognitive aspects of learning mathematics. A change in curriculum which brings affect into the forefront of undergraduate mathematics education could assuage the cognitive dissonance experienced by instructors and students alike who know or believe the process of learning and doing mathematics to be an emotional experience.

What if we actually included this deep emotional component of mathematics in undergraduate education? What form could this take? Couldn't we organize our teaching so that most students do actually experience such "aha" moments, and as frequently as possible?[9]

---

[9] To emphasize the creative process of doing mathematics, one simple idea (which both authors experimented with) is to include in a given course a series of assignments or questions focused on "The process of doing mathematics". Those may include identifying relationships between emotions and activity through emotional reflection and metacognition. Such assignments could be based on the framework of Mason, Burton, and Stacey's (1982) book "Thinking Mathematically", Chapter 2: Phases of Work.

An anarchist perspective of the undergraduate mathematics curriculum recognizes the need to constantly and meaningfully question the status quo to assess whether or not our current actions agree with our values. If the education community values recognizing and addressing emotion in the mathematics classroom, then an anarchist approach deems it necessary to build this into the curriculum.

By going through the creative process of mathematics in an undergraduate course, instead of – or in addition to – teaching a textbook manicured version, instructors can encourage students to try ideas. An inauthentic experience will not do, and some instructors may at first require recommendations for guiding students not only through the mathematical process, but also through its emotional connections. Furthermore, this is a process which takes time – and this is an idea that may put off those who adhere to a set of learning outcomes. Yet the experience of an authentic mathematical process is within reach of most undergraduate students. Most ideas will fail, and students will experience frustration. Then one will work, with its accompanying joy.

We argue that a byproduct of this emphasis on the emotional connection would be the development of resilience, perseverance, and the ability to stay focused and motivated. Those skills are essential for mathematics; in fact, many have argued that the skills required to be successful in mathematics and in endurance sports are very similar. This often comes as a surprise to non-mathematicians, as the inherent struggle of the mathematical process is purposely hidden and kept behind doors. By making it an integral part of our courses, we would help students build resilience in a supported environment, which is something that they could all benefit from.

In fact, pushing the analogy with endurance sports further, physical training also involves struggle to achieve a (usually self-directed) goal. Lots of sweating and sore muscles are involved in the process, often with coaching, but it's apparently an enjoyable activity. Why do people do this? Presumably they identify it as valuable – it makes them feel good, makes them more able to do the physical things they want to do. In undergraduate mathematics, could students embark on a similar journey? Undergraduate curriculum would need to provide direction and opportunities for affect to take a more prominent role in the learning process, and to explicitly recognize affect in the set of values of the mathematics education community.

Of course, if an instructor already holds these values, they may choose to make changes in their individual courses which elevate the emotional experiences involved in learning and doing mathematics. One small change could be simply to ask students to keep a journal of their emotional experiences in the class. This could be private or public, personal or anonymous, individual, or community-oriented, as agreed upon by the students and instructor. It may be that some students do not want to do this, and that is also fine. It is also likely that many undergraduate students will be pleasantly surprised for their emotions to be recognized and requested in a mathematics class. Changes to the curriculum can be large or small according to an anarchist perspective, so long as the changes continue to take place with the intention to align with values.

### *3.5. Creativity in Mathematics*

*We see mathematics as a creative endeavour, yet we mostly teach algorithmic procedures. We believe in the "throw ideas at the wall and see what sticks" creative process, yet we teach content that has been manicured to make it as clean and structurally sound as possible.*

Mathematicians tend to agree on an underlying force guiding their work. André Weil's research was art, discovery for Jaques Hadamard was driven by beauty, and for Poincare it was intuition. Contemporaries such as Claire Voisin and John Mason describe beauty and creativity as fundamental to the process of doing mathematics. Augustus de Morgan said that "[t]he moving power of mathematical invention is not reasoning but imagination" (de Morgan, 1866).

The creative process is messy and non-linear. What we read is the final product, not the process. This dichotomy is unsurprisingly akin to that of the artistic process. Indeed, a finished piece of artistry does not display the years required to build expertise and the messy process of production, just as the fancy restaurant does not leave its doors open to the kitchen.

Students however do not yet have this level of understanding, as they are barely starting to get familiar with mathematical thinking. We present them with the manicured version of mathematics; the polished, well structured statements that form the core of our modern textbooks. We show them established algorithms, and teach them how to apply them. How can we expect them to understand that this is not how mathematics is done? From what they see, they probably expect that to do mathematics is to come up with clear, concise, beautiful final solutions right from the start, which couldn't be further from the truth.

What if we instead taught mathematics like it is done? What if we exposed students to the mathematical process, and went through the process of throwing ideas at the wall with them? Creativity is an essential mathematical skill, and it is a skill that can be learned. Simple creative mathematical acts are described in mathematics education research, for example in the research that has grown from Kilpatrick's (1987) article on problem posing (see a review of the literature in Cai & Rott, 2024). Instead of presenting students with the final manicured version and asking them to understand and apply it, what if instructors went through the whole creative process with them, and have them formulate the final manicured version themselves? It is well-established that inquiry-based pedagogies equip students with a healthy range of skills for global citizenship (see, for example, the introduction to Archer-Kuhn, et al. (2020). It is also well-known that these pedagogies take more time and thus cover less content. What are the detriments to this way of learning, and would the benefits outweigh these?

> The mathematics curriculum is bursting at the seams. It has been characterized as a "mile wide and inch deep" [281]. Many mathematics courses, like calculus, have a rigid, scripted curriculum that leaves little room for creative thought or deep engagement. As a result, many mathematics professors spend all their time on "covering content" and spend little time to support the social aspects of learning. Yet, building relationships, community, and creating collective access for your students is the foundation of disrupting hierarchies. (Reinholz, 2023, p. 43)

As practical example, a culture of creativity can be developed in the classroom by asking solution space questions such as, "is this true always, sometimes, never?", generalization questions such as "what happens if we replace n with n+1?", or questions which are related but change the problem at hand to create a new, related problem. Some of these are explored in the literature on problem posing. For example, if the assigned problem is counting the number of graphs on four vertices, students can extend this question by asking "how many graphs are there on 5 vertices?" If the assigned problem asks for the possible solutions for a system of two equations in three unknowns, instructors can develop a mindset in the classroom in which students will extend the problem by asking about three equations in four unknowns. Instructors might choose to add the simple question, "change one small thing about this problem to create a new problem" to each assignment as a way to bring creative acts to the forefront of their classes.

## 4. Conclusion

We propose anarchism as a philosophical framework to help resolve the perceived cognitive dissonances that plague our undergraduate mathematics curriculum. Contemporary anarchism offers a perspective that focuses on diversity, pluralism, cooperation, autonomy and experimentation. Many mathematics practices, including the open-ended pursuit of unbiased truth, the celebration of diverse collaboration, and the shared responsibility of validating written work, are parallel to anarchist practices. It thus makes sense to reconsider undergraduate mathematics education using an anarchist framework.

Reflecting on values of diversity, open-endedness, and non-dominance, we see that mathematics education has a lot of progress to make. The focus of course and curriculum design on learning outcomes determined by a centralized few is inherently hierarchical and emphasizes the final product of mathematical thought over the process of learning and building mathematical skills. This is further reinforced by standard assessment practices which in turn marginalize the already disadvantaged and quell any diversity of thought or practice. It is true that the mathematics classroom is a uniquely difficult environment because of the tension between expertise and authority, but much progress could be made to move away from the hierarchical dominance that is reinforced in the classroom.

Mathematics education should serve multiple stakeholders. This means an invigorated focus on skills – not just basic cognitive skills but meta-cognitive skills, social and emotional skills, and practical and physical skills. Furthermore, skills should be integrated with knowledge, attitudes, and values. Although this is a perspective held by many, the cognitive dissonance between these values and the current mathematics curriculum spells out dissatisfaction among many students, faculty, employers, and other stakeholders.

It is also the authors' position that the human-ness of mathematics is ready to be brought into the curriculum in a substantial way. A focus on learning the practice of mathematics through skills, emotions, non-hierarchical teaching and curricular practices and most importantly a focus on processes rather than products can move us toward a more human and holistic mathematics education. The anarchist tenets of diversity and open-endedness are in line with the values of mathematicians and offer the undergraduate curriculum a path forward.

Ultimately, anarchism centers around an alignment between thoughts and actions. There is no anarchist blueprint of what an anarchist undergraduate curriculum should look like, and there should not be. Instead, in this perspective piece we use anarchism as a framework to provide food for thought for individual instructors' practice. We suggest both incremental and large-scale changes to the curriculum, moving from outcomes to process, from content to skills, and from hierarchical to horizontal teaching structures. As progress is made, we expect experiences to lead to new insights, new questions, and new paths forward.

We point you toward further reading, including a very recent and important book on hierarchy in mathematics (Reinholz, 2023) as well as persistent ideas about the psychology of learning mathematics (Skemp, 1971) in the hopes that you can enact your own anarchist progress toward improving mathematics education.

## References

Archer-Kuhn, B., Wiedeman, D., & Chalifoux, J. (2020). Student Engagement and Deep Learning in Higher Education: Reflections on Inquiry-Based Learning on Our Group Study Program Course in the UK. *Journal of Higher Education Outreach and Engagement*, *24*(2), 107-122.

Breitenbach, A., & Rizza, D. (2018). Introduction to special issue: Aesthetics in mathematics. *Philosophia Mathematica*, *26*(2), 153–160. https://doi.org/10.1093/philmat/nkx019

Cai, J., & Rott, B. (2024). On understanding mathematical problem-posing processes. *ZDM–Mathematics Education*, *56*(1), 61-71.

Cantor, G. (1883). Uber unendliche, lineare Punktmannigfaltigkeiten. *Mathematische Annalen 21*, 545-591.

Chandrasekhar S. (1987). Truth and Beauty. University of Chicago Press.

Chomsky, N. (1987). The Soviet Union versus Socialism. In D. I. Roussopoulos (Ed.), *The radical papers*, Black Rose Books, 49-51.

Chronaki, A. & Yolcu, A. (2021) Mathematics for "citizenship" and its "other" in a "global" world: critical issues on mathematics education, globalisation and local communities, Research in Mathematics Education, 23:3, 241-247.

de Morgan, A. (1866). Obituary of Sir William Rowan Hamilton. *The Gentleman's magazine*, vol. 220, p. 132. Online Version: https://babel.hathitrust.org/cgi/pt?id=mdp.39015027525990&seq=176/.

Elikai, F., & Schuhmann, P. W. (2010). An examination of the impact of grading policies on students' achievement. *Issues in Accounting Education, 25*(4), 677-693.

Freire, P. (2005). *Pedagogy of the oppressed*. Continuum.

Goldman, E. (1923). *My disillusionment in Russia*. New York: Doubleday, Page & Company. Online Version: https://www.marxists.org/reference/archive/goldman/works/1920s/disillusionment/index.htm .

Gordon, U. (2008). *Anarchy alive!: Anti-authoritarian politics from practice to theory*. Pluto Press.

Gu, T.E. & Bishop, M. (2022, Oct 14). Francis Su's "Mathematics on Human Flourishing": How Math Builds Virtues and Makes Us Human. The Phillipian. https://phillipian.net/2022/10/14/francis-sus-mathematics-on-human-flourishing-how-math-builds-virtues-and-makes-us-human/

Kilpatrick, J. (1987). Problem formulating: Where do good problems come from? In A. H. Schoenfeld (Ed.), *Cognitive science and mathematics education.* Hillsdale, 123-147.

Klymchuk, S., Sangwin, C. (2020). A new type of questions for teaching and assessing critical thinking in mathematics, Proceedings of the 12[th] International Conference of EDULEARN, 2348-2355.

Lockhart, P. (2002). *A mathematician's lament*. https://maa.org/sites/default/files/pdf/devlin/LockhartsLament.pdf .

Mason, J., Burton, L., & Stacey, K. (1982). *Thinking mathematically*. Addison Wesley.

McLeod, D. B. (1992). Research on affect in mathematics education: A reconceptualization. In D. A. Grouws (Ed.), *Handbook of research on mathematics teaching and learning: A project of the National Council of Teachers of Mathematics*. Macmillan Publishing Co, Inc, 575-596.

Nichols, T. (2017). *The death of expertise: The campaign against established knowledge and why it matters*. Oxford University Press.

Niss, M., Bruder, R., Planas, N., Turner, R., & Villa-Ochoa, J.A. (2017). In G. Kaiser (Ed.), Proceedings of the 13th International Congress on Mathematical Education, ICME-13 Monographs.

OECD Global Competencies (2018). *Preparing our youth for an inclusive and sustainable world: The OECD PISA global competence framework*. https://www.oecd.org/education/Global-competency-for-an-inclusive-world.pdf

Poincaré H. (1930). Science and Method. Maitland F. trans. London: Thomas Nelson.

Reinholz, D. (2023). *Equitable and engaging mathematics teaching: A guide to disrupting hierarchies in the classroom*. The Mathematical Association of America.

Remillard, J. T., & Heck, D. J. (2014). Conceptualizing the curriculum enactment process in mathematics education. ZDM Mathematics Education, 46(5), 705–718.

Rocker, R. (1989). *Anarcho-syndicalism: Theory and practice*. Pluto Press.

Skemp, R. (1971). *The psychology of learning mathematics*. Penguin, Harmondsworth.

Skovsmose, O., & Säljö, R. (2008). Learning mathematics through inquiry. *Nordic Studies in Mathematics Education, 13*(3), 31–52.

Society for Industrial and Applied Mathematics (SIAM) (2014). *Math modeling: Getting started and getting solutions*. https://live-siam-m3-challenge.pantheonsite.io/what-is-math-modeling/modeling-handbooks/

Society for Industrial and Applied Mathematics (SIAM) (2018). *Math modeling: Computing and communicating*. https://live-siam-m3-challenge.pantheonsite.io/what-is-math-modeling/modeling-handbooks/

Stommel, J. (2023). *Undoing the grade: Why we grade, and how to stop*. Hybrid Pedagogy.

Su, F. (2020). *Mathematics for human flourishing*. Yale University Press.